\documentclass[12pt]{article}

\usepackage{amsmath}
\usepackage{amssymb}

\newcommand{\llg}{{\rm lg}}

\newcommand{\eop}{\bigstar}  

\newcommand{\complB}{<\!\!\!\circ}
\newcommand{\card}[1]{{\vert #1 \vert} }

\newcommand{\forces}{\Vdash}

\newcommand{\crit}{{\rm crit}}

\newcommand{\Dom}{{\rm Dom}}

\newcommand{\PPr}{{\rm Pr}}
\newcommand{\Gr}{{\rm Gr}}

\newcommand{\Rang}{{\rm Rang}}

\newcommand{\Sup}{{\rm Sup}}

\newcommand{\cf}{{\rm cf}}

\newenvironment{Proof}{\noindent{\bf Proof.}}{\par\bigskip} 

\newtheorem{THEOREM}{Theorem}[section]

\newtheorem{Conclusion}[THEOREM]{Conclusion}

\newtheorem{Hypothesis}[THEOREM]{Hypothesis}

\newtheorem{LEMMA}[THEOREM]{Lemma}

\newtheorem{Main Theorem}[THEOREM]{Main Theorem}
\newenvironment{main Theorem}{\begin{Main Theorem}} 
{\end{Main Theorem}}

\newtheorem{Theorem}[THEOREM]{Theorem}

\newtheorem{Definition}[THEOREM]{Definition}

\newtheorem{Conventions}[THEOREM]{Conventions}

\newtheorem{Main Definition}[THEOREM]{Main Definition}
\newenvironment{main definition}{\begin{Main Definition}}
{\end{Main Definition}}

\newtheorem{Lemma}[THEOREM]{Lemma}

\newtheorem{Notation}[THEOREM]{Notation}

\newtheorem{Convention}[THEOREM]{Convention}

\newtheorem{Note}[THEOREM]{Note}

\newtheorem{Observation}[THEOREM]{Observation}

\newtheorem{Remark}[THEOREM]{Remark}

\newtheorem{Main Fact}[THEOREM]{Main Fact}
\newenvironment{main Fact}{\begin{Main Fact}}{\end{Main Fact}}

\newtheorem{Fact}[THEOREM]{Fact}

\newtheorem{Subfact}[THEOREM]{Subfact}

\newtheorem{Claim}[THEOREM]{Claim}

\newtheorem{Main Claim}[THEOREM]{Main Claim}
\newenvironment{main claim}{\begin{Main Claim}}{\end{Main Claim}}

\newtheorem{Crucial Claim}[THEOREM]{Crucial Claim}
\newenvironment{crucial claim}{\begin{Crucial Claim}}{\end{Crucial Claim}}

\newtheorem{Corrolary}[THEOREM]{Corrolary}

\newtheorem{Subclaim}[THEOREM]{Subclaim}

\newtheorem{Corollary}[THEOREM]{Corollary}

\newtheorem{Example}[THEOREM]{Example}

\newtheorem{Proposition}[THEOREM]{Proposition}

\newtheorem{Discussion}[THEOREM]{Discussion}

\newenvironment{Proof of the Subfact}
{\noindent{\bf Proof of the Subfact.}}{\par\bigskip}

\newenvironment{Proof of the Theorem}
{\noindent{\bf Proof of the Theorem.}}{\par\bigskip}

\newenvironment{Proof of the Conclusion}
{\noindent{\bf Proof of the Conclusion.}}{\par\bigskip}

\newenvironment{Proof of the Observation}
{\noindent{\bf Proof of the Observation.}}{\par\bigskip}

\newenvironment{Proof of the Fact}
{\noindent{\bf Proof of the Fact.}}{\par\bigskip}

\newenvironment{Proof of the Lemma}
{\noindent{\bf Proof of the Lemma.}}{\par\bigskip}

\newenvironment{Proof of the Claim}
{\noindent{\bf Proof of the Claim.}}{\par\bigskip}

\newenvironment{Proof of the Subclaim}
{\noindent{\bf Proof of the Subclaim.}}{\par\medskip}

\newenvironment{Proof of the Main Claim}
{\noindent{\bf Proof of the Main Claim.}}{\par\bigskip}

\newenvironment{Proof of the Crucial Claim}
{\noindent{\bf Proof of the Crucial Claim.}}{\par\bigskip}

\catcode`\@=11
\def\@begintheorem#1#2{\rm \trivlist \item[\hskip \labelsep{\bf #1\ #2.}]}
\def\@opargbegintheorem#1#2#3{\rm \trivlist
      \item[\hskip \labelsep{\bf #1\ #2\ (#3).}]}
\catcode`\@=12

\newcommand{\Bbf}{\Bbb}

\newcommand{\into}{\rightarrow}

\newcommand{\rest}{\upharpoonright}  

\newcommand{\deq}{\buildrel{\rm def}\over =}

\newcommand{\DD}{{\cal D}}
\newcommand{\EE}{{\cal E}}

\newcommand{\HH}{{\cal H}}
\newcommand{\II}{{\cal I}}
\newcommand{\JJ}{{\cal J}}
\newcommand{\KK}{{\cal K}}
\newcommand{\PP}{{\cal P}}

\newcount\skewfactor
\def\mathunderaccent#1#2 {\let\theaccent#1\skewfactor#2
\mathpalette\putaccentunder}
\def\putaccentunder#1#2{\oalign{$#1#2$\crcr\hidewidth
\vbox to.2ex{\hbox{$#1\skew\skewfactor\theaccent{}$}\vss}\hidewidth}}
\def\name{\mathunderaccent\tilde-3 }

\begin{document}

\title{Universal graphs at the successor of a singular cardinal}

\author{Mirna D\v zamonja\\
School of Mathematics\\
University of East Anglia\\
Norwich, NR4 7TJ, UK\\
\scriptsize{M.Dzamonja@uea.ac.uk}\\
\scriptsize{http://www.mth.uea.ac.uk/people/md.html}\\
and\\
Saharon Shelah\\
Institute of Mathematics\\
Hebrew University of Jerusalem\\
91904 Givat Ram, Israel\\
and\\
Rutgers University\\
New Brunswick, NJ, USA\\
\scriptsize{shelah@sunset.huji.ac.il}\\
\scriptsize{http://www.math.rutgers.edu/$\thicksim$shelarch}
}

\maketitle

\begin{abstract}

The paper is concerned with the existence of a universal graph at the
successor of a strong limit singular $\mu$ of cofinality $\aleph_0$.
Starting from the assumption of the existence of a supercompact
cardinal, a model is built in
which for some such $\mu$
there are $\mu^{++}$ graphs on $\mu^+$ that taken jointly are universal
for the graphs on $\mu^+$, while $2^{\mu^+}>>\mu^{++}$.

The paper also addresses
the
general problem of obtaining a framework
for consistency results at the successor of a singular strong limit
starting from the assumption that a
supercompact cardinal $\kappa$ exists. The
result on the existence of universal graphs is obtained as
a specific application of
a more general method.
{\footnote{This publication is denoted [DjSh 659] in Saharon Shelah's list of
publications. Both authors thank the United States-Israel Binational
Science Foundation for their support, and the NSF for their grant
NSF-DMS97-04477. Mirna D\v zamonja would like to acknowledge
that a major part of this research was done while she was a Van Vleck
Visiting Assistant Professor at the University of Wisconsin-Madison.

Keywords: successor of singular, iterated forcing, universal graph.

AMS 2000 Classification: 03E35, 03E55, 03E75.
}}
\end{abstract}

\baselineskip=16pt
\binoppenalty=10000
\relpenalty=10000
\raggedbottom

\section{Introduction}

The question of the existence of a universal graph
of certain cardinality and with certain properties has been the subject of
much research in mathematics (\cite{FuKo}, \cite{Kj}, \cite{KoSh},
\cite{random},
\cite{Sh 175a}, \cite{Sh 500}). By universality we mean here that 
every other graph of the same size embeds into the universal
graph. In the presence of $GCH$ it follows
from the classical results in model theory (\cite{ChKe}) that such a
graph exists at every uncountable cardinality, and it is well known
that the random graph (\cite{random}) is universal for 
countable graphs (although the situation is not so simple
when certain requirements on the graphs are imposed, 
see \cite{KoSh}). When the assumption of $GCH$ is dropped, it becomes
much harder to construct universal objects, and it is in fact
usually rather easy to obtain negative consistency results
by adding Cohen subsets to the universe (see \cite{KjSh 409}
for a discussion of this). For some classes of graphs
there are no universal objects as soon as $GCH$
fails sufficiently (\cite{Kj}), while 
for others there can exist consistently a small family of the class
that acts jointly as a universal object for the class at the
given cardinality (\cite{Sh 457}, \cite{DjSh 614}). 
Much of what is known in the absence of
$GCH$ is known about successors of regular cardinals
(\cite{Sh 457}, \cite{DjSh 614}). In \cite{Sh 175a} there is a positive
consistency result concerning the existence of a universal graph
at the successor of singular $\mu$ where $\mu$ is not a strong limit.
In this paper we address the issue
of the existence of a universal graph at the successor of a singular
strong limit and obtain a positive consistency results regarding the
existence of a small family of such graphs that act jointly as
universal for the graphs of the same size.

In addressing this specific problem,
the paper also offers a step towards the solution of a more
general problem of doing iterated forcing in connection with the
successor of a singular. This is the case because the result about universal
graphs is obtained as an application of a more general method. The 
method relies on an
iteration of $(<\kappa)$-directed-closed $\theta\ge\kappa^+$-cc forcing,
followed by the Prikry forcing for a normal ultrafilter $\DD$
built by the iteration. The cardinal
$\kappa$ here is supercompact in the ground model. The idea is that
the Prikry forcing for $\DD$ can be controlled by the iteration,
as $\DD$ is being built in the process
as the union
of an increasing sequence of normal filters that appear during
the iteration. Apart from building $\DD$, the iteration 
also takes care of the particular application it is aimed
at by predicting the $\DD$-names of the relevant
objects and taking care of them (in our application,
these objects are graphs on $\kappa^+$).
The iteration is followed by the Prikry forcing
for $\DD$, so changing the cofinality of $\kappa$ to $\aleph_0$.
Before doing the iteration we prepare $\kappa$ by rendering
its supercompactness indestructible by $(<\kappa)$-directed-closed
forcing through the use of Laver's diamond (\cite{Laver}). Not only
do we the use the indestructibility of $\kappa$, but 
Laver's diamond itself plays a crucial role in the definition
of the iteration. We note that the result has an unusual feature in
which the iteration is not constructed directly, but the existence
of such an iteration is proved and used.

Some of
the ideas connected to the forcing scheme discussed in this
paper were pursued by
A. Mekler and S. Shelah in \cite{MkSh}, and by M. Gitik and
S. Shelah in \cite{GiSh 597}, both in turn relaying on M. Magidor's independence
proof for $SCH$ at $\beth_\omega$ \cite{Magidor 1},
\cite{Magidor 2} and Laver's 
indestructibility method, \cite{Laver}. 
In \cite{MkSh}\S3 the idea of guessing Prikry names of an object
after the final collapse is present, while \cite{GiSh 597}
considers densities of box topologies, and for the particular forcing used
there presents a scheme similar to the one we use
(although the iteration is different). The latter paper also
reduced the strength of a large cardinal needed for the iteration
to a hyper-measurable. The difference between \cite{GiSh 597}
and our results is that the individual
forcing used in \cite{GiSh 597} is basically Cohen forcing, 
while our interest here is to 
give a general axiomatic framework under which the scheme can be applied for
many types of forcing notions.

The investigation of the consistent existence of universal
objects also has relevance in model theory. The idea here is to
classify theories in model theory by the size of their universality
spectrum, and much research has been done to confirm that this
classification is interesting from the model-theoretic point of
view (\cite{GrSh 174}, \cite{KjSh 409}, \cite{Sh 500}, \cite{DjSh 614}).
The results here sound a word of caution to this programme.
Our construction builds $\mu^{++}$ graphs on $\mu^+$
that are universal
for the graphs on $\mu^+$, while $2^{\mu^+}>>\mu^{++}$ and $\mu$
is a strong limit singular of cofinality $\aleph_0$. 
In this model we naturally obtain club guessing on $S^{\mu^+}_{\aleph_0}$
for order type $\mu$, and this will prevent
the prototype of a stable unsuperstable theory
${\rm Th}{}(^\omega\omega, E_n)_{n<\omega}$ from having a small universal
family, see \cite{Sh 457},
\cite{KjSh 447}. Hence the universality spectrum at such $\mu^+$
classifies the prototype of a simple unstable theory
(the theory of a random graph), as less complicated than
the prototype of a stable unsuperstable theory, contrary to the expectation.
A possible conclusion is that one should
concentrate the investigation of the universality spectrum as
a dividing line for unstable theories only on the
case $\lambda^+$ with $\lambda=\lambda^{<\lambda}$, as the case
of the successor of a singular is too sensitive to the set theory involved.

There are several further questions that this paper brings to mind.
From the point of view of model theory it would be
interesting to determine which other
first order theories fit the scheme of this paper and
from the point of view of graph theory one would like
to improve the result on the existence of $\mu^{++}$ jointly universal
graphs to having just one universal graph. Set-theoretically, we
would like to be able to replace
$\mu$ an unspecified singular strong limit by $\mu=\beth_\omega$, as well as
to investigate
singulars of different cofinality than $\aleph_0$. 
We did not concentrate here on obtaining
the right
consistency strength for our results, suggesting another question
that may be addressed in the future work.

The paper is organised as follows. 
The major issue is to define the iteration used in the second step 
of the above scheme,
which is done in certain generality in \S \ref{frame}. We give there
a sufficient condition for a one step forcing to fit the general scheme,
so obtaining an axiomatic version of the method. 
In \S\ref{graphs} we give the application to the existence of $\mu^{++}$
universal graphs of size $\mu^+$ for $\mu$ the successor of a strong
limit singular of cofinality $\aleph_0$.

Most of our notation is entirely standard, with the possible
exception of

\begin{Notation} For $\alpha$ and ordinal and a regular cardinal
$\kappa<\alpha$, we let
\[
S^\alpha_\kappa\deq\{\beta<\alpha:\,\cf(\beta)=\kappa\}.
\]
\end{Notation}

\section{The general framework for forcing}\label{frame} 

\begin{Definition} Suppose that $\kappa$ is a strongly inaccessible
cardinal $>\aleph_0$.
A function $h:\,\kappa\into \HH(\kappa)$ is called {\em Laver's diamond
on} $\kappa$
iff
for every $x$ and $\lambda$,
there is an elementary
embedding ${\bf j}:V\into M$ with
\begin{description}
\item{(1)} $\crit({\bf j})=\kappa$ and ${\bf j}(\kappa)>\lambda$,
\item{(2)} ${}^{\lambda}M\subseteq M$,
\item{(3)} $({\bf j}(h))(\kappa)=x$.
\end{description}
\end{Definition}

\begin{Theorem}\label{ldiam}{Laver} (\cite{Laver}) Suppose that
$\kappa$ is a supercompact cardinal. \underline{Then} there is a Laver's diamond
on $\kappa$.
\end{Theorem}

\begin{Hypothesis}\label{hipa} We work in a universe $V$ that satisfies
\begin{description}
 \item{(1)}
$\kappa$ is a supercompact cardinal, $\theta=\cf(\theta)\ge\kappa^+$ and
$GCH$ holds at and above $\kappa$,
\item{(2)}
$\Upsilon^{\theta}=\Upsilon\,\,\&\,\,\chi=\Upsilon^+$ and
\item{(3)}
$h:\,\kappa\into \HH(\kappa)$ is a Laver's diamond.
\end{description}
\end{Hypothesis}

\begin{Remark}
It is well known that the consistency of the above hypothesis follows from the
consistency of the existence of a supercompact cardinal. We in fact only use
the
$\chi$-supercompactness of $\kappa$.
\end{Remark}

\begin{Definition}\label{R}{Laver} (\cite{Laver}) We define
\[
\bar{R}=\langle R^+_\alpha, \name{R_\beta}:\,\alpha
\le\kappa, \beta<\kappa\rangle,
\]
an iteration done with Easton supports,
and a 
strictly increasing sequence $\langle \lambda_\alpha:\,\alpha<\kappa\rangle$
of cardinals,
where $\name{R_\alpha}$ and $\lambda_\alpha$ are defined by induction
on $\alpha <\kappa$ as follows.

If
\begin{description}
\item{(1)} $h(\alpha)=(\name{P},\lambda)$, where $\lambda$ is a
cardinal and $\name{P}$
is a $R^+_\alpha$-name of $(<\alpha)$-directed-closed forcing, and
\item{(2)} $(\forall \beta<\alpha)\,[\lambda_\beta<\alpha]$,
\end{description}
we let $\name{R}_\alpha\deq\name{P}$ and $\lambda_\alpha\deq\lambda$.
Otherwise, let $\name{R}_\alpha\deq\{\emptyset\}$ and $\lambda_\alpha
\deq\sup_{\beta<\alpha}\lambda_\beta$.

The extension in $R^+_\alpha$ is defined by letting
\[
p\le q\iff [\Dom(p)\subseteq\Dom(q)\,\,\&\,\,(\forall i\in\Dom(q))
(q\rest i\forces``p(i)\le q(i)")],
\]
(where $p$ denotes the weaker condition).
\end{Definition}

\begin{Remark}\label{credits} 
The forcing $R^+_\kappa$ used in this section is Laver's forcing
from \cite{Laver} which makes the
supercompactness of $\kappa$ indestructible under any
$(<\kappa)$-directed-closed forcing.
\end{Remark}

\begin{Convention}
Definitions \ref{init} and \ref{real},
Claim \ref{closed} and Observation \ref{usual} 
take place in $V_1\deq V^{R^+_\kappa}$. Notice that $\kappa^+\le\cf(\theta)
=\theta<\chi$ still holds in $V_1$, as $\Rang(h)\subseteq\HH(\kappa)$, and
that $\kappa$ is still supercompact.
\end{Convention}

\begin{Definition}\label{init} 
We define the family $\KK_\theta$ as the family of all sequences
\[
\bar{Q}=\langle P_i,\name{Q}_i,\name{A}_i:\,i<i^\ast=\llg(\bar{Q})<\chi\rangle
\]
which satisfy

\begin{description}
\item{(1)} $P_i\subseteq\HH(\chi)$ (and each $P_i$ is a forcing notion,
which will follow from the rest of the definition),
\item{(2)} $\langle P_i:\,i<i^\ast\rangle$ is $\complB$-increasing and each
satisfies $\chi$-cc,
\item{(3)} $\name{Q}_i$ is a $P_i$-name of a member of $\HH(\chi)$
(hence of cardinality $\le\Upsilon$),
\item{(4)} If $\cf(i)\ge\kappa$, then $P_i=\bigcup_{j<i}P_j$,
\item{(5)} $\name{A}_i$ is a canonical $P_{i+1}$-name of a subset of $\kappa$,
\item{(6)} Letting $G_i$ be $P_i$-generic over
$V_1$, then in $V_1[G_i]$,
\begin{description}
\item{(a)}
${\rm NUF}\deq\{\DD:\,\DD\mbox{ a normal ultrafilter on }\kappa\}$,
\item{(b)} for every $\DD\in {\rm NUF}$ we are given a
$(<\kappa)$-directed-closed forcing notion $Q^i_\DD \in
\HH(\chi)^{V_1[G_i]}$ whose minimal element is denoted by
$\emptyset_{Q^i_{\DD}}$,
\end{description}
\item{(7)} With the notation of (6),
we have that $\name{Q}_i[G_i]$ is 
\[
\{\emptyset\}
\cup {\rm NUF}\cup\left\{\{\DD\}\times Q^i_\DD:\,\DD\in {\rm NUF}\right\}.
\]

\item{(8)} The order on $\name{Q}_i[G_i]$ is given by letting 
\[
x\le y\mbox{ iff }[
x=y\mbox{ or }x=\emptyset\mbox{ or }(x=\DD\in {\rm NUF} \,\,\&\,\,
y\in \{x\}\times Q^i_\DD)\mbox{ or }
\]
\[
x=(\DD,x^\ast), y=(\DD,y^\ast) \mbox{ for some }\DD\in {\rm NUF}
\mbox{ and }Q^i_\DD\models `` x^\ast\le y^\ast"],
\]
\item{(9)}
We have
\[
P_i\deq\left\{p:\,\begin{array}{ll}
{\rm (i)}\, p\mbox{ is a function with domain }\subseteq i,\\
{\rm (ii)}\, j\in \Dom(p)\implies p(j)\mbox{ is a canonical }P_j
\mbox{-name }\\
\quad\mbox{                      of a member of }
\name{Q}_j\\
{\rm (iii)}\, \card{S\Dom(p)
}<\kappa\mbox{ (see below)}
\end{array}
\right\},
\]
ordered by letting
\[
p\le q\iff [\Dom(p)\subseteq\Dom(q)\,\,\&\,\,(\forall i\in\Dom(q))
(q\rest i\forces p(i)\le q(i))],
\]
where
\end{description}

({\em Definition \ref{init} continues below})

\begin{Notation} 
\begin{description}
\item{(A)}
For $i< i^\ast$,
and $p\in P_i$, we let
\end{description}

\begin{equation*}
S\Dom(p)\deq\left\{j\in \Dom(p):\,
\neg \left[p\rest j \forces_{P_j}
\begin{split}
``p(j)\in
\{\emptyset\}\cup \name{{\rm NUF}}\\
\cup\{(\name{\DD},\emptyset^j_{\name{\DD}}):\,\DD\in \name{NUF}\}"\\
\end{split}
\right]\right\},
\end{equation*}

\begin{description}
\item{(B)} For $i< i^\ast$ and $p\in P_i$
we call $p$ {\em purely full} iff:

$S\Dom(p)=\emptyset$ and
for every $j<i$ we have
\[
p\rest j\forces_{P_j}``p(j)\in \name{{\rm NUF}}".
\]
\item{(C)}
Suppose that $i<i^\ast$ and $p\in P_i$ is purely full,
we define
\[
P_i/p\deq\{q\in P_i:\,q\ge p\,\,\&\,\,\mbox{each }q(j)\mbox{ has form
}(\name{\DD},x)\mbox{ for some }x\},
\]
with the order inherited from $P_i$.
\end{description}
\end{Notation}

({\em Definition \ref{init} continues:})

\begin{description}
\item{(10)} For every $i\le i^\ast$ and $p\in P_i$ which is purely full
we have that $P_i/p$ satisfies $\theta$-cc
and $P_i/p\in \HH(\chi)$.

\end{description}

\end{Definition}

\begin{Observation} (1) If $\bar{Q}\in \KK_\theta$ and $i<\llg(\bar{\theta})$,
\underline{then} $P_{i+1}=P_i\ast\name{Q}_i$.

{\noindent (2)} Assuming that $\langle P_j,\name{Q}_j,\name{A}_j:\,j<i\rangle
\in \KK_\theta$ and (5)-(10) above hold, we can see that 
$\langle P_j,\name{Q}_j,\name{A}_j:\,j\le i\rangle\in \KK_\theta$. Hence
$\KK_\theta$ can be alternatively defined by specifying when $\langle P_j,\name{Q}_j,\name{A}_j:\,j<i\rangle
\in \KK_\theta$ by induction on $i<\chi$.
\end{Observation}

\begin{Definition}\label{real}
\begin{description}
\item{(1)}
Let $\kappa^+\le\cf(\theta)=\theta<\chi$.
We define the family $\KK^+_\theta$ as the family of all sequences
\[
\bar{Q}=\langle P_i,\name{Q}_i,\name{A}_i:\,i<\chi\rangle
\]
such that
\[
i<\chi\implies \bar{Q}\rest i\in \KK_\theta.
\]
We let $P_\chi\deq\bigcup_{i<\chi}P_i$.

\item{(2)} Suppose that $\bar{Q}\in  \KK^+_\theta$ and
$\langle p_i:\,i<\chi\rangle$ 
with $p_i\in P_{\zeta_i}$ are purely full and increasing in $P_\chi$, where
$\zeta_i\deq\min\{\zeta:\,p_i\in P_{\zeta}\}$ (so if $i<j$ then
$p_i=p_j\rest \zeta_i$).
We define 
\[
P_\chi/\cup_{i<\chi} p_i\deq\left\{q\in P_\chi:\,
(\forall i<\chi)[q\rest \zeta_i \in P_{\zeta_i}/p_i]\right\},
\]
with the order inherited from $P_\chi$.
\end{description}
\end{Definition}

\begin{Claim}\label{closed}
\begin{description}
\item{(1)}
If $\bar{Q}\in \KK_\theta$, \underline{then} for all $i\le \lg(\bar{Q})$, we have
that $P_i$ is $(<\kappa)$-directed-closed.

\item{(2)} Similarly for $\bar{Q}\in \KK_\theta^+$.
\end{description}
\end{Claim}

\begin{Proof of the Claim} (1) Given 
a directed family $\{p_\alpha:\,\alpha<\alpha^\ast<\kappa\}$ of conditions in $P_i$. 
We shall define a common extension $p$ of this family. Let us first let
$\Dom(p)\deq\bigcup_{\alpha<\alpha^\ast} \Dom(p_\alpha)$. For $j\in
\Dom(p)$, we define $p(j)$ by induction on $j$. We work in $V_1^{P_j}$
and assume that $\{p_\alpha\rest j:\,\alpha<\alpha^\ast\}\subseteq G_{P_j}$.

If $j\notin \bigcup_{\alpha<\alpha^\ast} S\Dom(p_\alpha)$, then notice that there
is at most one $\DD\neq\emptyset$ such that for some (possibly more than
one)
$\alpha<\alpha^\ast$ we have $p_\alpha\rest j\forces``p_\alpha(j)=\name{\DD}"$,
as the family is directed. If there is such $\DD$, we let $p(j)
\deq\DD$, otherwise we let $p(j)=\emptyset$.

If $j\in \bigcup_{\alpha<\alpha^\ast} S\Dom(p_\alpha)$, similarly to the last
paragraph, we conclude that there is exactly one $\DD$ such that
\[
[\alpha<\alpha^\ast\,\,\&\,\,j\in S\Dom(p_\alpha)]\implies
p_\alpha\rest j\forces``p_\alpha(j)\in\{\name{\DD}\}\times\name{Q}^j_{\DD}".
\]
As $\name{Q}^j_{\DD}$ is forced to be $(<\kappa)$-directed-closed, we can
find in $V^{P_j}_1$ a condition $q$ such that $q\ge p_j(\alpha)$ for all
$\alpha<\alpha^\ast$ such that $j\in S\Dom(p_\alpha)$.
Let $p(j)\deq q$ for some such $q$.

{\noindent (2)} Follows from (1) as $\chi=\cf(\chi)>\kappa$.
$\eop_{\ref{closed}}$
\end{Proof of the Claim}

\begin{Observation}\label{usual} Suppose that $\bar{Q}\in \KK^+_\theta$,
$i<j<\chi$ and $p\in P_i$, $q\in P_j$
are purely full, while $p\le q$. \underline{Then}

\begin{description}
\item{(1)} $\Dom(p)\subseteq \Dom(q)$ and $\alpha\in \Dom(p)\implies
p(\alpha)=q(\alpha)$.
\item{(2)}
Suppose that $r\in P_i/p$.

\underline{Then} defining $r+q\in P_j$ by letting $\Dom(r+q)=
\Dom(q)$ and letting for $\alpha\in \Dom(r)$
\[
(r+q)(\alpha)\deq\left\{
\begin{array}{ll}
r(\alpha) &\mbox{ if }\alpha\in \Dom(p)\\
(q(\alpha), \emptyset_{\name{Q}^\alpha_{q(\alpha)}}) &\mbox{ otherwise,}
\end{array}
\right.
\]
we obtain a condition in $P_j/q$.

\item{(3)} For $r_1,r_2\in P_i/p$ we have that
\begin{description}
\item{$(\alpha)$}
$r_1$ and $r_2$ are
incompatible in $P_i/p$ iff $r_1+q$ and $r_2+q$ are incompatible
in $P_j/q$,
\item{($\beta)$} $r_1\le_{P_i/p}r_2\iff r_1+q\le_{P_j/q}
r_2+q$.
\end{description}

\item{(4)} $P_i/p\complB_f P_j/q$ where $f(r)\deq r+q$.

\item{(5)}
Suppose
that the sequence
$\bar{p}=\langle p_i:\,i<\chi\rangle$ satisfies 
that each $p_i\in P_{\zeta_i}$ is purely full,
and the sequence $\bar{p}$ is increasing in $P_\chi$,
where
\[
\zeta_i\deq\min\{\zeta:\,p_i\in P_\zeta\}=\xi_i+1 >i
\]
and $\langle \zeta_i:\,i<\chi\rangle$ is strictly increasing.
\underline{Then}  $P^\ast=P_\chi/\cup_{i<\chi}p_i$
is isomorphic to the limit of a
$(<\kappa)$-supported
iteration of $(<\kappa)$-directed-closed $\theta$-cc
forcing,
namely
\[
P^\ast\approx
\lim\langle
P_{\xi_i}/(p_i\rest \xi_i),\name{Q}^{\xi_i}_{p_i(\xi_i)}:\,i<\chi\rangle,
\]
with the complete embeddings $f_{p_i,p_j}:\,P_{\xi_i}/(p_i\rest \xi_i)
\into P_{\xi_j}/(p_j\rest \xi_j)$ as in (4) above.

\item{(6)} For every $r\in P_\chi$, there is $q\ge r$ with
$S\Dom(q)=S\Dom(r)$ and $p$ purely full in some $P_i$, such that $q\in P_i/p$.
\end{description}
\end{Observation}

\begin{Convention}\label{fs} Since $f_{p_i,p_j}$ are usually clear form the
context we simplify the notation by not mentioning these functions explicitly.
\end{Convention}

\begin{Claim}\label{names} Suppose that $\bar{Q}\in \KK^+_\theta$ and
$\name{t}$ is a $P_\chi$-name of an ordinal, while $p\in P_\chi$ is purely full.

\underline{Then} for some $j<\chi$ and $q$ we have
$p\le q\in P_j$, and $q$ is purely full, and above $q$
we have that $\name{t}$ is a $P_j$-name
(i.e $\name{t}$ is a $P_j/q$-name).
\end{Claim}

\begin{Proof of the Claim} Given $p\in P_\chi$ purely full, and suppose that
the conclusion fails. Let $i<\chi$ be such that $p\in P_i$.  We 
shall choose by induction on $\zeta<\theta$
ordinals $i_\zeta$ and $\gamma_\zeta$ and condition $r_\zeta$
such that
\begin{description}
\item{(i)}  $i_\zeta\in [i,\chi)$ and $\langle _\zeta:\,\zeta<\theta\rangle$ is
increasing continuous, \item{(ii)} $p_\zeta\in P_{i_\zeta}$ is purely full, with
$p_0=p$,
\item{(iii)} $p_{\zeta}\le r_\zeta$ with
$r_\zeta\forces_{P\chi}``\name{t}=\gamma_\zeta"$,
\item{(iv)} $\gamma_\zeta\notin\{\gamma_\xi:\,\xi<\zeta\}$, so in particular
$r_\zeta$ is incompatible with every $r_\xi$ for $\xi<\zeta$,
\item{(v)} $p_\zeta\deq\cup_{\xi<\zeta}p_\xi$ for $\zeta$ a limit.
\item{(vi)}  
$r_\zeta\in P_{i_{\zeta+1}}/p_{\zeta+1}$. \end{description}
We now explain how to do this induction.

Given $p_\zeta$ and $i_\zeta$.
Since we are assuming that $\name{t}$
is not a $P_{i_\zeta}$-name
above $p_\zeta$, it must be possible to find
$r_\zeta$ and $\gamma_\zeta$ as required. 
Having chosen $r_\zeta$,
(by extending $r_\zeta$ if necessary),
we can choose $p_{\zeta+1}$ as required
in item (vi) above, see Observation \ref{usual}(6).
This determines $i_{\zeta+1}$. Note that $i_{\zeta+1}<\chi$ as
$P_\chi\deq\bigcup_{j<\chi} P_j$.

However, completing the induction we arrive at a contradiction,
as letting
$p^\ast
\deq\cup_{\zeta<\theta}p_\zeta$ we obtain
a purely full condition. Hence $P\deq P_{\sup_{\zeta<\theta}i_\zeta}/p^\ast$
has $\theta$-cc, but $\{ r_\zeta+p^\ast:\,\zeta<\theta
\}$ forms
a set of $\theta$ pairwise incompatible conditions in
$P$.
$\eop_{\ref{names}}$
\end{Proof of the Claim}

\begin{Convention} Now we go back to $V$, i.e. the Main Claim \ref{point}
takes place in $V$.
\end{Convention}

\begin{Main Claim}\label{point} Suppose
\begin{description} 
\item{($\alpha$)} $\bar{\name{Q}}=\langle \name{P}_i,\name{Q}_i,\name{A}_i:\,
i<\chi\rangle$
is an $R^+_\kappa$-name for a member of $\name{\KK}^+_\theta$,

\item{$(\beta$)} ${\bf j}:\,V\into M$ is an elementary embedding
such that ${}^{\Upsilon}M\subseteq M$, $\crit({\bf j})=\kappa$,
$\chi<{\bf j}(\kappa)$ and
\[
({\bf j}(h))(\kappa)=(\name{P}_\chi,\chi)
\]
{\scriptsize{(such a choice is possible by the definition of Laver's diamond)
.}} \end{description}

Considering ${\bf j}(\langle R^+_\alpha,\name{R}_\alpha:\,
\alpha<\kappa\rangle)$ in $M$, by its definition we see that
\[
{\bf j}(\langle R^+_\alpha,\name{R}_\alpha:\,\alpha<\kappa\rangle)
=\langle R^+_\alpha,\name{R}_\alpha:\,\alpha<{\bf j}(\kappa)\rangle
\]
and $\name{R}_\kappa=\name{P}_\chi$.
Hence ${\bf j}(R^+_\kappa)=R^+_\kappa\ast\name{P}_\chi\ast\name{R}^\ast$
for some $R^+_\kappa\ast\name{P}_\chi$-name $\name{R}^\ast\in M$
for a forcing notion, which is forced to be $\chi^+$-closed.

We also let
\[
\bar{\name{Q}}'=\langle \name{P}'_i,\name{Q}'_i,\name{A}'_i:\,
i<{\bf j}(\chi)\rangle\deq {\bf j}(\langle \name{P}_i,\name{Q}_i,\name{A}_i:\,
i<\chi\rangle).
\]

\underline{Then} in $V^{R^+_\kappa}$, 
the following
holds:
we can find $\bar{\alpha}=
\langle\alpha_i:\,i<\chi\rangle$,
$\bar{p}^\ast=\langle p_i^\ast:\,i<\chi\rangle$ and
$\bar{q}^\ast=\langle q_i^\ast=({}^1 q_i,{}^2 q_i):\,i<\chi\rangle$ such that

\begin{description}
\item{(a)} $\langle \alpha_i:\,i<\chi\rangle$ is strictly increasing
continuous and each $\alpha_i<\chi$,

\item{(b)} $p^\ast_i\in P_{\alpha_i+1}$ is purely full,

\item{(c)} $\bar{p}^\ast$ is increasing in $P_\chi$,

\item{(d)} For every $i<\chi$, we have $\bar{q}^\ast\rest i
\in M^{R^+_\kappa}$, and in $M^{R^+_\kappa}$ we have
\[
(p^\ast_i, {}^1q_i,{}^2q_i)
\in P_\chi\ast \name{R}^\ast
\ast \name{P}'_{{\bf j}(\alpha_i+1)},
\]
while $(p^\ast_i, {}^1q_i)\in P_\chi\ast
 \name{R}^\ast$,

\item{(e)} In $M^{R^+_\kappa}$ we have that for $\gamma<\chi$
\[
\langle (p^\ast_i,{}^1q_i,{}^2q_i)
:\,i<\gamma\rangle\mbox{ is increasing in }
P_\chi\ast\name{R}^\ast
\ast
{\name{P}}'_{\sup_{i<\gamma}{\bf j}(\alpha_i+1)},
\]

\item{(f)} In $M^{R^+_\kappa}$,
it is forced by $(p^\ast_{i+1}, {}^1q_{i+1})$ that
${}^2q_{i+1}$ is an upper bound to
\[
\{{\bf j}(r):\,r\in 
\name{G}_{P_{\alpha_i}\ast
\name{Q}^{\alpha_i+1}_{p^\ast_{i}(\alpha_i)}}\}, \]

\item{(g)} If $\name{B}$ is an $R^+_\kappa$-name of a
$P_{\name{\alpha}_i+1}$-name of a subset of $\kappa$,
\underline{then}
for some $R^+_\kappa\ast\name{P}_{\chi}$-name
$\name{{\bf t}}_{\name{B}}$ for a truth value (i.e. an ordinal $\in
\{0,1\}$):
\begin{description}

\item{(1)} In $V$ we have that
$(\emptyset_{R^+_\kappa},p^\ast_{i+1})$ forces $\name{{\bf t}}_{\name{B}}$
to be a $P_{\name{\alpha}_{i+1}+1}/p^\ast_{i+1}$-name,

\item{(2)}
$M\models [(\emptyset_{R^+_\kappa},p^\ast_{i+1},q^\ast_{i+1})\forces``
\kappa\in{\bf j}(\name{B})\mbox{ iff }\name{{\bf t}}_{\name{B}}=1"].$
\end{description}

\item{(i)} In $M^{R^+_\kappa}$, either
\[
(p_{i+1}^\ast,q^\ast_{i+1})
\forces``\kappa\in {\bf j}(\name{A}_{\alpha_i})",
\]
or:
\[
p^\ast_i
\forces_{
P_\chi}
\begin{array}{c}
``\mbox{ there is no }q=({}^1q,{}^2q)
\ge_{\name{R}^\ast\ast \name{P}'_{{\bf j}(\alpha_i)+1}} q^\ast_i
\mbox{ with }\\
{}^1q\forces_{\name{R}^\ast}`` {}^2q({\bf{j}}(\alpha_i))\ge_{\name{P}'_{{\bf
j}(\alpha_i)+1}} \{{\bf j}(r):\,r\in
{\name{G}}_{P_{\alpha_i}\ast\name{Q}^{\alpha_i}_{p_i^\ast(\alpha_i)}}\}\\
\mbox{ and }
\kappa\in {\bf j}(\name{A}_{\alpha_i})"".
\end{array}
\]

[Note that ${\bf j}(\name{A}_{\alpha_i})$ is a $P'_{{\bf
j}(\alpha_i)+1}$-name for a subset of ${\bf j}(\kappa)$.]

\item{(j)} If $\cf(i)\ge \theta$,
\underline{then} in $V^{R^+_\kappa\ast\name{P}_{\name{\alpha}_i}}$ we have
$p^\ast_i(\alpha_i)\in NUF$ and specifically
\end{description}
\begin{equation*}
p^\ast_{i}(\alpha_i)=\left\{
\name{B}[G_{P_{\alpha_i}}]:\,
\begin{split}
\name{B}\mbox{ is a }P_{\alpha_i}/(p^\ast_i\rest\alpha_i)\mbox{-name
for a subset of }\kappa\\
\quad\mbox{ and }\name{{\bf t}}_{\name{B}}[G_{P_{\alpha_i}}]=1\\
\end{split}
\right\}.
\end{equation*}

\end{Main Claim}

\begin{Remark} In fact,
to accommodate various applications, we might want to
weaken item (j) of the Main Claim \ref{point},
say to apply only to stationary many $i\in S^\chi_{\ge\theta}$. The same proof
would work, but as we do not need this at present, we shall not go into this
generality.
\end{Remark}

\begin{Proof of the Main Claim}
Consider $\langle R^+_i,\name{R}_i:\,\kappa<i<{\bf j}(\kappa)\rangle$
over $R^+_\kappa\ast\name{P}_\chi$ in $M$.
By the inductive definition of $R_i$
(which is preserved by {\bf j}), for $\kappa<i<\chi^+$, we have
that $\name{R}_i$ is a name for the trivial forcing. For $\chi^+<i
<{\bf j}(\kappa)$, we have that $\name{R}_i$ is a name for a
$(<\chi)$-directed-closed forcing in $M$, so in $V$ as well, as
${}^{<\chi}M\subseteq M$. Similarly we conclude that
$\name{P}'_{{\bf j}(\zeta)}$ names a $(<\chi)$-closed forcing notion,
for all $\zeta<\chi$.
This observation will be used repeatedly and in particular
will enable us to use the master condition idea in the induction
below. In particular, we can conclude that $\name{R}_\chi$ is
$(<\chi)$-complete. By the choice of {\bf j},
\[
\forces_{{\bf j}(R^+_\kappa)}``\mbox{each }P'_i/p\mbox{ is
}(<\chi)\mbox{-closed for }p\in P'_i\mbox{ purely full}."
\]
Also notice that in the induction below, we have
that in $V_1$, the cardinality of $P_{\alpha_i+1}/p^\ast_i$
is $\le\Upsilon$ and $P_{\alpha_i+1}/p_i^\ast$ satisfies $\theta$-cc,
so in $V_1^{P_{\alpha_i}}$ we have $2^\kappa\le\Upsilon$.

Now we choose $(\alpha_i,p^\ast_i,q^\ast_i)$ in $M^{R^+_\kappa}$
by induction on $i$. We start with $\alpha_0=0$, $p_0^\ast\in P_1$
any purely full condition, and $q^\ast_0=\emptyset$.

\underline{Choice of $p^\ast_{i+1}, q^\ast_{i+1}$ and $\alpha_{i+1}$}.

Given $p^\ast_i$ and $\alpha_i$ in $V^{R^+_\kappa}$. We have that
\[
p^\ast_i\rest\alpha_i\forces_{P_{\alpha_i}}
``\card{\name{Q}^{\alpha_i}_{p^\ast_i(\alpha_i)}}\le\Upsilon
\,\,\&\,\,
\name{G}_{\name{Q}^{\alpha_i}_{p^\ast_i(\alpha_i)}}
\subseteq \name{P}_{\alpha_i+1}/p^\ast_i."
\]
Hence in $M$, letting $\name{X}_i\deq\{{\bf j}(r):\,r\in 
\name{G}_{P_{\alpha_i}\ast\name{Q}^{\alpha_i}_{p^\ast_i(\alpha_i)}}\}$ we have
\begin{equation*}
(\emptyset_{R^+_{{\bf j}(\kappa)}},{\bf j}(p^\ast_i\rest\alpha_i))
 \forces_{P'_{{\bf j}(\alpha_i)}}
\begin{split}
``\name{X}_i
\subseteq \name{P}'_{{\bf j}(\alpha_i)+1}
\mbox{ is directed}\\
\mbox{ and above}\\
{\bf j}(p_i(\alpha_i))\mbox{ has size}\le\Upsilon."\\
\end{split}
\end{equation*}
In $V_1$, we have that the forcing $P_{\alpha_i+1}/p^\ast_i$
is a $\theta$-cc forcing notion of size $\le\Upsilon$,
hence there are $\le\Upsilon^\theta\cdot\Upsilon=\Upsilon$
canonical $P_{\alpha_i+1}/p^\ast_i$-names for a subset of
$\kappa$. Let us enumerate them as
$\langle \name{B}_\zeta^{i+1}:\,\zeta<\zeta^\ast(i+1)
\le\Upsilon\rangle$, with $\name{B}^{i+1}_0=\name{A}_i$.
By induction on $\zeta\le\zeta^\ast(i+1)$ we choose purely full $p^{i+1}_\zeta$
increasing continuous with $\zeta$,
$q^{i+1}_\zeta=({}^1 q^{i+1}_\zeta,
{}^2 q^{i+1}_\zeta)$ increasing with $\zeta$,
$\alpha_\zeta^{i+1}$ increasing with $\zeta$
and $\name{{\bf t}}_{\name{B}_\zeta^{i+1}}$ as follows.

Let $p^{i+1}_0\deq p^\ast_i$, $\alpha^{i+1}_0\deq \alpha_i$ and
$q^{i+1}_0\deq q^\ast_i$. 

Coming to $\zeta+1$, let $G$ be a $R^+_\kappa\ast\name{P}_\chi$
generic such that $(\emptyset_{R^+_\kappa}, p^{i+1}_\zeta)\in G$
and let $H$ be a {\bf j}($R^+_\kappa\ast\name{P}_\chi)$
generic over $M$ so that
$\{{\bf j}(r):\,r\in G\}\subseteq H$. This can be achieved by the
familiar argument using the fact that ${\bf j}(R^+_\kappa\ast\name{P}_\chi)$
is $(<{\bf j}(\kappa))$-directed-closed, while
$G$ is $(<\kappa)$-directed and has size $\le\chi<{\bf j}(\kappa)$. In
particular, ${\bf j}$ lifts to an embedding of $V[G]\into M[H]$.
By the fact that ${\bf j}(R^+_\kappa\ast\name{P}_\chi)=
R^+_\kappa\ast\name{P}_\chi\ast\name{R}^\ast\name{P}'_{{\bf j}(\chi)}$,
we can write $M[H]=M[H_0][H_1]$ where $H_0$ is $R^+_\kappa\ast\name{P}_\chi$
generic over $M$.
In $M[H_0]$ we ask
``the $\zeta$-question":

Is it true that
there is no 
\[
q=({}^1 q, {}^2 q)
\ge_{\name{R}^\ast\ast \name{P}'_{{\bf j}(\alpha_\zeta^{i+1})+1}}
\{q^{i+1}_\xi:\,\xi\le\zeta\}
\]
with 
\[
{}^1 q\forces_{\name{R}^\ast}``{}^2 q\ge\name{X}_i
\,\,\&\,\,\kappa\in{\bf j}(\name{B}^{i+1}_\zeta)"
\,\,\&\,\,{}^2 q\in\name{P}'_{{\bf j}(\alpha^{i+1}_\zeta)+1}/
({\bf j}(p^{i+1}_\zeta)\rest{\bf j}(\alpha^{i+1}_\zeta)+1)?
\]
If the answer is positive, in $M$ we define 
$\name{{\bf t}}_{\name{B}^{i+1}_\zeta}
\deq 0$ (hence a $R^+_\kappa\ast\name{P}_\chi$-name for a truth
value),
and
\[
q^{i+1}_{\zeta+1}=({}^1 q^{i+1}_{\zeta+1},{}^2 q^{i+1}_{\zeta+1})
\]
to be any $R^+_\kappa\ast\name{P}_\chi$-name for a condition in
$\name{R}^\ast\ast \name{P}'_{{\bf j}(\chi)}$ such that
\[
(\emptyset_{R^+_\kappa}, p^{i+1}_\zeta)\forces``{}^1 q^{i+1}_{\zeta+1}
\ge_{\name{R}^\ast}{}^1 q^{i+1}_{\xi}"
\]
for every $\xi\le\zeta$,
and 
\[
(\emptyset_{R^+_\kappa}, p^{i+1}_\zeta,{}^1 q^{i+1}_{\zeta+1})
\forces_{\name{P}'_{{\bf j}(\alpha^{i+1}_\zeta)+1}}
``{}^2 q^{i+1}_{\zeta+1}\ge\cup\{{}^2
q_\xi^{i+1}:\,\xi\le\zeta\}
\,\,\&\,\,{}^2 q^{i+1}_{\zeta+1}\ge \name{X}_i".
\]
The choice of ${}^1 q^{i+1}_{\zeta +1}$ is possible by the
induction hypothesis and the fact that 
\[
\forces_{R^+_\kappa\ast \name{P}_\chi}``\name{R}^\ast\mbox{ is }
(<\chi)\mbox{-directed-closed}".
\]
Let us verify that the choice of ${}^2 q^{i+1}_{\zeta+1}$ is possible.
Working in $M$ we have that
$(\emptyset_{R^+_\kappa}, p^{i+1}_\zeta,{}^1 q^{i+1}_{\zeta+1})$
forces $\name{X}_i$ over 
$(\emptyset_{R^+_\kappa}, p^{i+1}_\zeta,{}^1 q^{i+1}_{\zeta+1})$
to be a $(<\kappa)$-directed subset of $\name{P}'_{{\bf j}(\chi)}$ of size
$\le\chi$. Hence
if $\zeta=0$ we can choose ${}^2 q^{i+1}$
to be forced to be above $\name{X}_i$. We can similarly
choose ${}^2 q^{i+1}_{\zeta+1}$ for $\zeta>0$.

If the answer to the $\zeta$ question
is negative, we let $\name{{\bf t}}_{\name{B}_\zeta^{i+1}}\deq 1$
and choose
$q^{i+1}_{\zeta+1}=({}^1 q^{i+1}_{\zeta+1},{}^2 q^{i+1}_{\zeta+1})$
in $M$ exemplifying the negative answer.

At any rate, $\name{{\bf t}}_{\name{B}_\zeta^{i+1}}$ is a $R^+_\kappa
\ast\name{P}_\chi$-name for an ordinal. By Claim \ref{name},
in $V^{R^+_\kappa}$
there is $\alpha^{i+1}_{\zeta+1}\ge \alpha^{i+1}_\zeta$ and purely full
$p^{i+1}_{\zeta+1}\ge p^{i+1}_\zeta$ with
$p^{i+1}_{\zeta+1}\in P_{\alpha^{i+1}_{\zeta+1}}$ such that
$\name{{\bf t}}_{\name{B}_\zeta^{i+1}}$ is a $P_{\alpha^{i+1}_{\zeta+1}}/
p^{i+1}_{\zeta+1}$-name.

For $\zeta$ limit, let $\alpha^{i+1}_\zeta\deq\sup_{\xi<\zeta}\alpha^{i+1}_\xi$,
$p^{i+1}_\zeta\deq\cup_{\xi<\zeta}p^{i+1}_\xi$, and
$q^{i+1}_\zeta$ not defined.

At the end, we let
$\alpha_{i+1}\deq\sup_{\zeta<\zeta^\ast(i+1)}\alpha^{i+1}_{\zeta+1}$
and $p^\ast_{i+1}$ any purely full condition in $P_{\alpha_{i+1}+1}$ with
$p^\ast_{i+1}\ge \bigcup_{\zeta<\zeta^\ast(i+1)}p^{i+1}_\zeta$,
and $q^\ast_{i+1}$ such that
\[
(\emptyset_{R^+_\kappa}, p^\ast_{i+1})\forces``q^\ast_{i+1}
\ge_{\name{R}^\ast\ast\name{P}'_{{\bf j}(\alpha_{i+1})+1}}
\{q^{i+1}_\zeta:\,\zeta\le\zeta^\ast(i+1)\}".
\]

\underline{Choice of $p^\ast_i, q^\ast_i$ and $\alpha_i$ for $i<\chi$ limit}.
We let $\alpha_i\deq\sup_{j<i}\alpha_j$ and choose $p^\ast_i
\in P_{\alpha_i+1}$ purely full so that $p^\ast_i\ge \cup_{j\le i}p^\ast_j$,
and if $\cf(i)\ge\theta$, \underline{then}
\begin{equation*}
p^\ast_i(\alpha_i)\deq\left\{
\begin{split}
\name{B}[G_{P_{\alpha_i}}]:\,
&\name{B}\mbox{ is a }P_{\alpha_i}/(p^\ast_i\rest\alpha_i)\mbox{-name for a
subset of }
\kappa\\
&\mbox{ and }\name{{\bf t}}_{B}[G_{P_{\alpha_i}}]=1\\
\end{split}
\right\}.
\end{equation*}
It follows by the construction of Laver's diamond and standard arguments about
elementary embeddings and master conditions that
 \[
p^\ast_{i}\rest\alpha_i\forces_{P_{\alpha_i}}``
p^\ast_i(\alpha_i)\in\name{{\rm NUF}}".
\]
Then we can choose $q^\ast_i$ so that
$(\emptyset_{R^+_\kappa}, p^\ast_i, q^\ast_i)\ge (\emptyset_{R^+_\kappa},
p^\ast_j, q^\ast_j)$ for all $j<i$
and $q^\ast_i\ge\{{\bf j}(r):\,r\in G_{P_{\alpha_i}}\}$, which is again
possible by
the observation at the beginning of the proof.
$\eop_{\ref{point}}$
\end{Proof of the Main Claim}

\begin{Conclusion}\label{glavna}
In $V_1$, if $\bar{Q}\in \KK^+_\theta$
and $\langle p^\ast_i:\,i<\chi\rangle$ is as guaranteed by Main Claim
\ref{point},
letting $\name{\DD}_i\deq p^\ast_{i}(\alpha_i)$, it follows by Observation
\ref{usual}(5) that 
\[
\bar{P}=\langle P_{\alpha_i}/(p^\ast_i\rest
\alpha_i),\name{Q}^{\alpha_i}_{\name{\DD}_i}:\,
i<\chi\rangle
\]
is an iteration with $(<\kappa)$-supports
of 
$(<\kappa)$-directed-closed $\theta$-cc forcing. In addition,
there is a club $C$ of $\chi$ with the property
that in $V^{P}_1$
\[
\langle \DD_i:\,i\in C\,\,\&\,\,\cf(i)\ge\theta
\rangle
\]
is an increasing sequence of normal filters
over $\kappa$, with
\[
[i\in C\,\,\&\,\,\cf(i)\ge\theta]
\implies P_{\alpha_i}/(p^\ast_i\rest\alpha_i)\forces``
\name{\DD}_i\mbox{ is an ultrafilter over }\kappa".
\]
If $\delta<\chi$ satisfies $\cf(\delta)>\kappa$ then $\cup_{i<\delta} p^\ast_i$ forces
over $P_{\alpha_\delta}$ that $\bigcup_{i<\delta}
\name{\DD}_i$ us an ultrafilter over $\kappa$ which is generated by
$\cf(\delta)$ sets.
\end{Conclusion}

\begin{Definition}\label{fitted} (In $V^{R^+_\kappa}$) Given
$\bar{Q}=\langle P_i,\name{Q}_i,\name{A}_i:\,i<\chi\rangle
\in \KK^+_\theta$.

We say that $\bar{Q}$ is {\em fitted} iff there is a continuous
increasing sequence $\langle \alpha_i:\,i<\chi\rangle$ of ordinals
$<\chi$, and a sequence $\langle p^\ast_i:\,i<\chi\rangle$ of purely full
conditions  with $p^\ast_i\in P_{\alpha_i+1}$, such that
letting $\name{\DD}_i\deq p^\ast_{i(\alpha_i)}$,
\[
\langle P_{\alpha_i+1}/(p^\ast_i\rest
\alpha_i),\name{Q}^{\alpha_i}_{\name{\DD}_i}:\,
i<\chi\rangle
\]
is an iteration with $(<\kappa)$-supports of 
$(<\kappa)$-directed-closed $\theta$-cc forcing,
\underline{and}
\[
\cf(i)\ge\theta\implies\forces_{P_{\alpha_i+1}/(p^\ast_i\rest\alpha_i)}
``\name{A}_i\in\name{\DD}_i".
\]
\end{Definition}

\begin{Crucial Claim}\label{reflect}
(In $V^{R^+_\kappa}$) The following is a sufficient
condition for $\name{\bar{Q}}\in\KK^+_\theta$ to be fitted:

There is a definition ${\bf R}$ such that:
\begin{description}
\item{(1)} for every
forcing ${\Bbf P}$ with $\card{{\Bbf P}}\le\Upsilon$ in $V^{\Bbf P}$
and a ${\Bbf P}$-name $\name{\DD}$ of a normal ultrafilter on $\kappa$
we have that ${\bf R}[{\Bbf P},\name{\DD}]$ is a ${\Bbf P}$-name of a forcing
notion of cardinality $\le\Upsilon$,
\item{(2)} for every
purely full $p\in P_\chi$ and $i\in \Dom(p)$, we have
that
\[
p\rest i \forces``p(i)=\name{Q}^i_{p(i)}={\bf R}[P_i/(p\rest
i),p(i)]",
\]
\item{(3)} there is a definition $\name{f}$ that for
every
forcing ${\Bbf P}$ with $\card{{\Bbf P}}\le\Upsilon$ in $V^{\Bbf P}$
and a ${\Bbf P}$-name $\name{\DD}$ of a normal ultrafilter on $\kappa$
gives a ${\Bbf P}$-name of a
function ${\name{f}}_{[{\Bbf P},\name{\DD}]}:\,{\bf R}[{\Bbf P},\name{\DD}]
\into\name{\DD}$ such that
for every purely full $p\in P_\chi$ and $i\in \Dom(p)$
it is forced by $p\rest i$ that:

``for every inaccessible $\kappa'<\kappa$ and every ${\bf g}$
a  $(<\kappa')$-directed family of conditions in ${\bf R}[P/(p\rest i),p(i)]$
of size $<\kappa$, such that 
\[
r\in {\bf g}\implies\kappa'\in {\name{f}}_{[P/(p\rest i), p(i)]}(r),
\]
\underline{there is} $q\ge {\bf g}$ such that
$q\forces\kappa'\in \name{A}_i$." 

\end{description}
\end{Crucial Claim}

\begin{Remark} The condition in Claim \ref{reflect} is sufficient
for the present application in \S2.
It may be weakened if needed for some future
application. Really, the condition to use instead
of it is that in item (i) of Main Claim \ref{point}, for all $i$ of cofinality
$<\theta$, we are ``in the good case", i.e. the first case of item (i).
However, we
wish to have a criterion which can be used without the knowledge of the proof
of the Main Claim \ref{point}, and the condition in Claim \ref{reflect} is one
such criterion.



\end{Remark}

\begin{Proof of the Crucial Claim}
By Conclusion \ref{glavna} it suffices to show that under the assumptions of
this Claim,
in the proof of Main Claim \ref{point} we can choose
$\langle \alpha_i:\,i<\chi\rangle$,
$\langle p^\ast_i:\,i<\chi\rangle$ and
$\langle q^\ast_i:\,i<\chi\rangle$ so that
for every $i$ with $\cf(i)\ge \theta$, the answer to
``the 1st question" in the choice of $q^{i+1}_1$ is negative.
The proof is by induction on such $i$. We use the notation of Main Claim
\ref{point}.

Given $i$ with $\cf(i)\ge\theta$. Hence we have
\begin{equation*}
p^\ast_i(\alpha_i)=\left\{
\name{B}[G_i]:\,
\begin{split}
\name{B}\mbox{ a }P_{\alpha_i}/(p^\ast_i\rest\alpha_i)
\mbox{-name for a}\\
\mbox{subset of }\kappa\mbox{ and }
\name{{\bf t}}_{\name{B}}=1\\
\end{split}
\right\}\deq\name{\DD}_i.
\end{equation*}
In $M$ we have

\begin{equation*}
(\emptyset_{R^+_\kappa}, p^\ast_i, q^\ast_i)
\forces
\left(
\begin{split}
\{{\bf j}(r)({\bf j}(\alpha_i)):\,{\bf j}(r)\in \name{X}_i\}
\mbox{ is }(<\kappa)-\mbox{directed of size }<\\
{\bf j}(\kappa),\kappa\mbox{ is inaccessible and }
(\forall r)[\kappa\in {\bf j}({\name{f}}_{[P/(p^\ast\rest\alpha_i),
p^\ast_i(\alpha_i)]}(r))]\\
\end{split}
\right).
\end{equation*}
(The last statement is true by the definition of $\name{\DD}_i$ and
$\name{t}_{\name{B}}$,
no matter what $f_{[P/p^\ast\rest\alpha_i,
p^\ast_i(\alpha_i)]}(r)$ is forced to be.)

By the assumption (3) and elementarity, applying ${\bf j}$ we have that the
answer to the ``1st question" is negative.
$\eop_{\ref{reflect}}$
\end{Proof of the Crucial Claim}

\begin{Definition}\label{nova} (In $V^{R^+_\kappa}$)
Given $\theta=\cf(\theta)\in (\kappa,\chi)$.
We define $\KK^\ast_\theta$ in the same way as $\KK^+_\theta$,
but with a freedom of choice for $Q_0$. Namely, to obtain
the definition of $\KK^\ast_\theta$ from that of
$\KK^+_\theta$, we
\begin{description}
\item{(A)} In item (6) of Definition \ref{init}, require
$i>0$,
\item{(B)} We let $Q_0$ be any $(<\kappa)$-directed-closed
$\theta$-cc forcing notion in
$\HH(\chi)$.
\end{description}
\end{Definition}

\begin{Claim}\label{OK} (In $V^{R^+_\kappa}$)
Main Claim \ref{point}, Conclusion \ref{glavna},
Definition \ref{fitted}
and Claim \ref{reflect} hold with
$\KK^+_\theta$ replaced by $\KK^\ast_\theta$.
\end{Claim}

\begin{Proof of the Claim} As in $V^{R^+_\kappa\ast\name{Q}_0}$,
$\kappa$ is still indestructibly supercompact and $\Upsilon^{\theta}
=\Upsilon$.
$\eop_{\ref{OK}}$
\end{Proof of the Claim}

\begin{Discussion}\label{success}
\begin{description}
\item{(1)}
In the present application, we need to make sure that
cardinals are not collapsed, so we 
have $\theta=\kappa^+$ and is $Q_{\DD}$ chosen to have a
strong version of $\kappa^+$-cc which is preserved by iterations
with $(<\kappa)$-supports. 
\item{(2)}
Clearly, Claim \ref{reflect} remains true if we replace the
word
``inaccessible" by e.g ``strongly inaccessible", ``weakly compact",
``measurable".
\item{(3)} As we shall see in section \ref{graphs}, the point of dealing with a
fitted member of $\KK^+_\theta$  is to be able to control the Prikry names in
the forcing that will be performed after the iteration extracted from
$\KK^+_\theta$, namely the Prikry forcing over $\cup_{i<\delta}\DD_i$ for some
$\delta$.
The point of $\name{A}_i$ is to give us a control of this ultrafilter
in the appropriate universe. With this in mind,
we could use Claim \ref{reflect} to represent our results in the axiomatic
form, and there is also an equivalent game-theoretic representation. As it is
not entirely clear that Claim \ref{reflect} is the best sufficient condition
for fittedness, we have decided not to formulate any axioms here.
\end{description}
\end{Discussion}

\section{Universal graphs}\label{graphs}

\begin{Theorem}\label{graph}
Assume that it is consistent that a supercompact cardinal $\kappa$
exists, and let $\Upsilon$ and $\chi$ be such that $\Upsilon^+=\chi$ and
$\Upsilon^{\kappa+}=\Upsilon$. 

\underline{Then} it is consistent to have a singular
strong limit cardinal
$\mu$ of cofinality $\omega$ with $2^{\mu^+}=\chi>\mu^{++}$,
on which there are $\mu^{++}$ graphs of size $\mu^+$ which are
universal for the graphs of size $\mu^+$.
\end{Theorem}

\begin{Proof}
We start with a universe $V$ in which $\mu$,
$\Upsilon$ and $\chi$ satisfy Hypothesis 
\ref{hipa}, with $\mu$ in place of $\kappa$ and $\theta=\mu^+$. Let
$R^+_\kappa$ be the forcing described in Definition \ref{R}. We work in
$V^{R^+_\kappa}$, which we start calling $V$ from this point on. As we shall
not use $R^+_\kappa$ any more, we free the notation $\name{R}_\alpha$ to be
used with a different meaning in this section.

\begin{Definition}\label{Q0}
Let $Q_0$
be the Cohen forcing which
makes $2^{\mu^+}=\Upsilon$ by adding $\Upsilon$ distinct $\mu^+$-branches
$\{\eta_\alpha:\,\alpha<\Upsilon\}$ to $({}^{\mu^+>}2)^{V}$ by conditions of
size $\le\mu$. Let $V_0
\deq V[G_{Q_{0}}]$.
\end{Definition}

\begin{Notation} If $\kappa$ is measurable and $\DD$ is a normal ultrafilter on
$\kappa$,  let $\PPr(\DD)$ denote the Prikry forcing for $\DD$.
\end{Notation}

\begin{Discussion} The idea of the proof is to embed ``$\DD$-named
graphs" into a universal graph. We use an iteration of forcing to achieve this.
As we intend to perform a Prikry forcing at the end of iteration, we need to 
control the names of graphs that appear \underline{after} the Prikry forcing,
so one
worry is that there would be too many names to take care of by the bookkeeping.
Luckily, we shall not be dealing with all such names, but only with those for
which we  are sure that they will actually be used at the end. This is achieved
by building the ultrafilter that will serve for the Prikry forcing, as the
union of filters that appear during the iteration. To this end, for every
relevant $\DD$ we also force
a set $\name{A}$ that  will in some sense be a ``diagonal
intersection" of $\DD$, so its membership in the intended ultrafilter will
guarantee that that ultrafilter contains $\DD$ as a subset.
\end{Discussion}

\begin{Definition}\label{defa} Suppose $V'\supseteq V_0$ 
is a universe in which $2^{\kappa^+}\le\Upsilon$, while $\kappa\le\mu$ is
measurable and $\DD$ is a normal ultrafilter over $\kappa$.
Working in $V'$,
we define a forcing notion $Q=Q_{\DD}\deq Q^{V'}_{\DD,\kappa}$, as follows.

Let $\bar{\name{M}}=\langle\name{M}_\alpha
=\langle \kappa^+, \name{R}_\alpha\rangle:\,\alpha<\Upsilon\rangle$ list
without repetitions
all canonical (in the usual sense) $\PPr(\DD)$-names for graphs on $\kappa^+$.
For definiteness we pick the first such list in the canonical well-order of
$\HH(\chi)$.
Elements of $Q$ are of the form
\[
p=\langle A^p, B^p, u^p,
\bar{f}^p=\langle f^p_\alpha:\,\alpha\in u^p\rangle\rangle,
\]
where
\begin{description}
\item{(i)} $A^p\in [\kappa]^{<\kappa}$,
\item{(ii)} $B^p\in \DD\cap \PP([\kappa\setminus(\Sup(A^p))])$,
\item{(iii)} $u^p\in [\Upsilon]^{<\kappa}$,
\item{(iv)} For $\alpha\in u^p$, we have that $f^p_\alpha$ is a partial
one-to-one function from $\kappa^+$ with $\card{\Dom(f^p_\alpha)}<\kappa$,
mapping $\zeta\in\Dom(f^p_\alpha)$
to an element of $\{\eta_\alpha\rest\zeta\}\times\kappa$,
\item{(v)} For $\alpha,\beta\in u^p$, for every $x', x'', y', y''$,
if 
\[
f^p_\alpha(x')=f^p_\beta(y')\neq f^p_\alpha(x'')=f^p_\beta(y''),
\]
\underline{then} for every $w\in [A^p]^{<\aleph_0}$
\[
\langle w,B^p\rangle
\forces_{\PPr(\DD)}``\name{M}_\alpha\models\name{R}_\alpha
(x',x'')\mbox{ iff }
\name{M}_\beta\models\name{R}_\beta
(y',y'')",
\]
\underline{and} for every $w\in [A^p]^{<\aleph_0}$ the 
condition $\langle w,B^p\rangle$
decides in the Prikry forcing for $\DD$ if
$\name{M}_\alpha\models\name{R}_\alpha (x',x'')$
\end{description}
We define the order on $Q$ by letting $p\le q$ (here $q$ is a stronger
condition) iff
\begin{description}
\item{(a)} $A^p$ is an initial segment of $A^q$,
\item{(b)} $A^q\setminus A^p\subseteq B^p$,
\item{(c)} $B^p\supseteq B^q$,
\item{(d)} $u^p\subseteq u^q$,
\item{(e)} For $\alpha\in u^p$, we have $f^p_\alpha\subseteq
f^q_\alpha$.
\end{description}
\end{Definition}

\begin{Claim}\label{properties} Suppose that $Q=Q_{\DD,\kappa}^{V'}$
is defined
as in Definition \ref{defa}. \underline{Then} in $V'$:
\begin{description} 
\item{(1)} $Q$ is a separative partial order.
\item{(2)} Suppose that $G$ is $Q$-generic over $V'$, and let
in $V'[G]$
\[
A^\ast\deq\bigcup \{A:\,(\exists B, u, \bar{f})[\langle A,B,u,\bar{f}
\rangle\in G]\}.
\] 
\underline{Then} $A^\ast\in [\kappa]^{\kappa}$ and $A^\ast
\subseteq^\ast B$ for every $B\in \DD$.
\item{(3)} For $\alpha<\Upsilon$ and $a\in \kappa^+$, the set
\[
\KK_{a,\alpha}\deq\{p\in Q:\,\alpha\in u^p\,\,\&\,\, a\in \Dom(f^p_\alpha)\}
\]
is dense open in $Q$.
\end{description}
\end{Claim}

\begin{Proof of the Claim} {\noindent (1)} Routine checking.
\begin{description}
\item{(2)} For $\alpha<\kappa$, the set
\[
\II_\alpha\deq\{p\in Q:\,(\exists\beta\ge\alpha)[\beta\in A^p]\}
\]
is dense open, hence $A^\ast\in [\kappa]^\kappa$. For $B\in \DD$ the set
\[
\JJ_B\deq\{p\in Q:\,B^p\subseteq B\}
\]
is dense open. If $p\in {\JJ}_B\cap G$, then for any $q\in G$
with $q\ge p$ we have $A^q\setminus A^p\subseteq B^p$.
Hence $A^\ast\setminus B\subseteq A^p$.

\item{(3)} Given $p\in Q$, clearly there is $q\ge p$ with $\alpha
\in u^q$. Without loss of generality $\alpha\in u^p$ and $a\notin
\Dom(f^p_\alpha)$. Applying
the Prikry Lemma, for every $b\in\Dom(f^p_\alpha)$ and 
$w\in [A^p]^{<\aleph_0}$, there is $B_{w,b}\subseteq B^p$ with
$B_{w,b}\in \DD$ and such that
\[
(w, B_{w,b})||_{\PPr(\DD)} ``\name{M}_\alpha\models b \name{R}_\alpha a".
\]
Choose $\gamma<\kappa$ such that $(\eta_\alpha\rest a,\gamma)
\notin \bigcup_{\beta\in u^p}\Rang(f^p_\beta)$, which is possible
as for every relevant $\beta$ we have $\card{\Dom(f^p_\beta)}<
\kappa$. Now we define $q$ by letting $A^q\deq A^p$, $B^q\deq
\bigcap\{B_{w,b}:\,w\in [A]^{<\aleph_0}\,\,\&\,\,b\in \Dom(f^p_\alpha)
\}$,
$u^q\deq u^p$ and
\[
f^q_\beta\deq\left\{
\begin{array}{lr}
f^p_\beta &\mbox{ if }\beta\neq\alpha\\
f^p_\alpha\cup\{(a,(\eta_\alpha\rest a,\gamma)\} &\mbox{ otherwise.}
\end{array}
\right.
\]
\end{description}
$\eop_{\ref{properties}}$
\end{Proof of the Claim}

\begin{Notation} Suppose that $Q$ is as in Claim \ref{properties}.
For $\alpha<\Upsilon$ let
\[
\name{f}_\alpha\deq\cup\{f^p_\alpha:\,\alpha\in u^p\,\,\&\,\,p\in\name{G}_Q\}.
\]
\end{Notation}

\begin{Definition} (Shelah, \cite{Sh 80}) Let $\lambda\ge\aleph_0$
be a cardinal. A forcing notion $P$ is said to be
{\em stationary} $\lambda^+$-{\em cc} iff for every $\langle
p_\alpha:\,\alpha<\lambda^+\rangle$ in $P$, there is a club
$C\subseteq\lambda^+$ and a regressive $f:\,\lambda^+
\into\lambda^+$ such that for all $\alpha,\beta\in C$,
\[
[\cf(\alpha)=\cf(\beta)=\lambda
\,\,\&\,\,f(\alpha)=f(\beta)]\implies p_\alpha, p_\beta
\mbox{ are compatible}.
\]
\end{Definition}

\begin{Theorem}{Shelah}\label{fact} (\cite{Sh 80})
Suppose that $\lambda^{<\lambda}=\lambda\ge\aleph_0$.
Iterations with $(<\lambda)$-support
of $(<\lambda)$-directed-closed stationary $\lambda^+$-cc
forcing, are $(<\lambda)$-directed-closed and satisfy stationary $\lambda^+$-cc.
\end{Theorem}

\begin{Claim}\label{cc} Suppose that $Q$ is as in Claim
\ref{properties}. \underline{Then} $Q$ is $(<\kappa)$-directed-closed
and satisfies stationary $\kappa^+$-cc.
\end{Claim}

\begin{Proof of the Claim} First suppose that $i^\ast<\kappa$
and $\{p_i:\,i<i^\ast\}$ is directed. For $i<i^\ast$
let $p_i\deq \langle A^i, B^i, u^i, \bar{f}^i\rangle$.
We define $A\deq\bigcup_{i<i^\ast}A^i$, $B\deq\bigcap_{i<i^\ast}
B^i$, $u\deq\cup_{i<i^\ast} u^i$, and for $\alpha\in u$ we let
$f_\alpha\deq\cup_{i<i^\ast}f^i_\alpha$. It is easily verified that
this defines a common upper bound of all $p_i$.

Hence $Q$ is $(<\kappa)$-directed-closed. Now we shall prove that
it is $\kappa^+$-stationary-cc. Let $\langle p_i:\,i<\kappa^+
\rangle$ be given, where
each $p_i=\langle A^i, B^i, u^i, \bar{f}^i
\rangle$. 

There is a stationary $S\subseteq S^{\kappa^+}_\kappa$
and $A^=\in [\kappa]^{<\kappa}$
and $\sigma,\tau<\kappa$ such that for all $i\in S$ we have $A^i=A^=$
and $\card{u^i}=\sigma$, and
$\card{\bigcup\{\Dom(f^i_\alpha):\,\alpha\in u^i\}}=\tau$.
For $i\in S$, let $\zeta_i\deq\sup\bigcup_{\alpha\in u^i}
\Dom(f^i_\alpha)$, hence $\cf(\zeta_i)<\kappa$. So
\[
E\deq\{j<\kappa^+:\,\cf(j)\ge\kappa\implies(\forall i<j)[\zeta_i<j]\}
\]
is a club of $\kappa^+$. Let $S_1\deq S\cap E$.

Let $\theta\deq\sigma+\tau+\card{A^=}$, so $\theta<\kappa$.
For $i\in S_1$ let $u^i\deq\{\alpha^i_s:\,s<\sigma\}$ be an increasing
enumeration. For every such $i$, we define a model $M_i$ with universe
$\kappa^+$, relations ${\bf R}^i_{w,s}$ for $w\in [A^=]^{<\aleph_0}$
and $s<\sigma$, and (partial) functions $g^i_s$
from $\kappa^+$ to $\kappa$, for $s<\sigma$.
This model is defined by letting
\[
(\zeta,\xi)\in {\bf R}^i_{w,s}
\mbox{ iff }[\zeta,\xi\in \Dom(f^i_{\alpha^i_s})
\,\,\&\,\, (w,B^i)\forces_{\PPr(\DD)}``\zeta\name{R}_{\alpha^i_s}
\xi"],
\]
and
\[
g^i_s(\zeta)=\gamma \mbox{ iff }f^i_{\alpha^i_s}(\zeta)=
(\eta_{\alpha^i_s},\gamma).
\]
Note that always
$\card{{\bf R}^i_{w,s}}\le\theta$ and $\card{g^i_s}
\le\theta$.

Now let $X$ be the set of all isomorphism types of models
with their universe an ordinal $<\kappa^+$ and $\le\theta$
relations and functions, each of cardinality $\le\theta$.
Hence $\card{X}=\kappa^+$, let $X\deq\{t_i:\,i<\kappa^+\}$.
Now note that there is a club $C$ of $\kappa^+$ such that
for every $j\in C\cap S^{\kappa^+}_\kappa$, types of all
models with universe $< j$ and $\le\theta$ relations of functions,
each of cardinality $\le\theta$, are enumerated in $X$ with
an index $<j$. Let $S_2\deq C\cap S_1$.

For $i\in S_2$, let $h(i)=l$ iff $t_l$ is the type of
$M_i\rest \sup (i\cap\bigcup_{\alpha\in u^i}\Dom(f_\alpha^i))+1$.
Hence $h$ is regressive on $S_2$, so there is a stationary
subset $S_3$ of $S_2$ such that $h$ is constant on $S_3$.

It is easily verified
that $p_i, p_j$ are compatible for every $i,j\in S_3$.
$\eop_{\ref{cc}}$
\end{Proof of the Claim}

\begin{Observation}\label{ufs}
Suppose that $\DD$
is a normal ultrafilter over $\kappa$
and $Q$ is a forcing notion such that
\[
\forces_Q``\DD\subseteq\name{\DD}'
\mbox{ and }\name{\DD}'\mbox{ is a normal ultrafilter
over }\kappa".
\]
\underline{Then} $\PPr(\DD)\complB_f
Q\ast\PPr(\name{\DD}')$, where $f$ is the
embedding given by 
\[
f((a, A))\deq (\emptyset_Q, (a,A)).
\]
\end{Observation}

\begin{Definition}\label{gr} Suppose that $Q$ is as in Claim \ref{properties},
while $Q\complB P$, and $\name{\DD}'$ is a $P$-name of a normal
ultrafilter over $\kappa$, extending $\DD\cup\{\name{A}^\ast\}$.
For $\alpha<\Upsilon$ we define $\name{Gr}^{{\DD}'}_\alpha$, intended to 
be a name 
for a graph on $\{\eta_\alpha\rest\zeta:\,\zeta<\kappa^+\}\times\kappa$
(see Claim \ref{name} below),
defined by letting for $y',y''\in \{\eta_\alpha\rest\zeta:\,\zeta<\kappa^+\}
\times\kappa$,
\[
\begin{array}{ccc}
y'\name{R} y''&\mbox{ iff for some }\langle p,\langle w, B^p\rangle\rangle\in
\name{G}\mbox{ with } \alpha\in u^p, p\in Q\mbox{ and }[w]\in [A^p]^{<\aleph_0}\\
&\mbox{ and some }x',x''\in\Dom(f^p_\alpha)\\
&\mbox{ we have } f^p_\alpha(x')=y'\mbox{ and }f^p_\alpha(x'')=y'',\\
&\mbox{ AND }\langle w, B^p\rangle\forces_{\PPr(\name{\DD})}
``\name{M}_\alpha\models\name{R}_\alpha(x',x'')".
\end{array}
\]
\end{Definition}


\begin{Claim}\label{name}  Suppose $Q$ is as in Claim \ref{properties},
while $Q\complB P$, and $\name{\DD}'$ is a $P$-name of a normal
ultrafilter over $\kappa$, extending $\DD\cup\{\name{A}^\ast\}$
(equivalently, $\name{A}^\ast\in\name{\DD}'$).

\underline{Then}
\[
\langle\emptyset, \langle\emptyset,\name{A}^\ast\rangle\rangle
\forces_{P\ast\PPr(\name{\DD}')}``\name{f}_\alpha \mbox{ is an embedding of }
\name{M}_\alpha  \mbox{ into } \name{Gr}^{\name{\DD}'}_\alpha.
\]
\end{Claim}

\begin{Proof of the Claim} Let $G$ be $P\ast\PPr(\name{\DD}')$-generic
with $\langle\emptyset, \langle\emptyset,\name{A}^\ast\rangle\rangle\in G$
and suppose that $x', x''$ are such that $M_\alpha\models R_\alpha(x', x'')$
in $V[G]$. Let $\langle p^+, \langle w,\name{A}'\rangle\rangle$ be a condition
in $G$ that
forces this. Without loss of generality, we have
\[
\langle p^+, \langle w,\name{A}'\rangle\rangle\ge
\langle\emptyset, \langle\emptyset,\name{A}^\ast\rangle\rangle.
\]
In particular,
$p^+\forces_P``w\in [\name{A}^\ast]^{<\aleph_0}"$. Considering $P$ as
$Q\ast P/Q$, let us write $\langle p^+, \langle w,\name{A}'\rangle\rangle$
as $\langle p, p', \langle w,\name{A}'\rangle\rangle$.
As $\name{A}^\ast$ is a $Q$-name, by extending $p^+$ if necessary, we may
assume that $A^p\supseteq w$, and then using the density of $\KK_{x',\alpha}$
and $\KK_{x'',\alpha}$, we may also assume that $\alpha\in u^p$ and $x',
x''\in\Dom(f^p_\alpha)$. By extending further, we may assume that $p^+\forces
``\name{A}'\subseteq B^p"$. Then $\langle p^+,\langle w, \name{A}'\rangle\rangle
\ge \langle p, \langle w, B^p\rangle\rangle$, hence the latter is in $G$.
Since $p\forces_P``\langle w, b^p\rangle ||_{\PPr(\DD)}\name{R}_\alpha(x',
x'')"$, it must be that $\langle w, B^p\rangle\forces_{\PPr(\DD)}
``M_\alpha\models R_\alpha(x', x'')"$. Hence in $V[G]$ we have that 
\[
y'=f_\alpha(x') R y''=f_\alpha(x'').
\]
On the other hand, suppose that in $V[G]$ we have 
$y'=f_\alpha(x') R y''=f_\alpha(x'')$ and let $\langle p, \langle w,
B^p\rangle\rangle$ exemplify this. In particular, $\langle w, B^p\rangle$
forces in $\PPr(\DD)$ that $``M_\alpha\models R_\alpha(x', x'')"$, and since 
$\langle p, \langle w,
B^p\rangle\rangle\in G$, we have that $R_\alpha(x', x'')$ holds in $V[G]$.

As it is easily seen that each $f_\alpha$ is forced to be 1-1 and total, this
finishes the proof.
$\eop_{\ref{name}}$ \end{Proof of the Claim}

\begin{Claim}\label{universal} Suppose that $Q$ and $\name{\DD}'$
are as in Claim \ref{name},
while $G$ is $Q$-generic over $V'$ and $2^\kappa=\kappa^+$ holds in $V'$.
Further suppose that $H$ is a $\PPr(\DD')$-generic filter over $V'[G]$ with $\langle \emptyset,
A^\ast\rangle\in H$.

\underline{Then} in $V'[G][H]$, there is a graph $Gr^\ast$
of size
$\kappa^+$ such that for every
$\PPr(\DD)$-generic filter $J$ over $V'$,
every graph of size $\kappa^+$ in $V'[J]$
is embedded into $Gr^\ast$.
\end{Claim}

\begin{Proof of the Claim} Define $Gr^\ast$ 
on
$\cup_{\alpha<\Upsilon}
\{\eta_\alpha\rest\zeta:\,\zeta<\kappa^+\}\times\kappa$,
hence $\card{Gr^\ast}=\kappa^+$, by our assumptions on $V'$.
We let
\[
Gr^\ast
\models``(\eta_\alpha\rest\zeta,i)R(\eta_\alpha\rest\xi,j)" 
\mbox{ iff }
Gr_\alpha^{\DD'}\models``(\eta_\alpha\rest\zeta,i)R(\eta_\alpha
\rest\xi,j)".
\]
Then $Gr^\ast$ is a well defined graph, as follows by
the definition of $Q$.

Given $M$ a graph on $\kappa^+$ in $V'[J]$, there is $\alpha$
such that $M=\name{M}_\alpha[G][J]$, hence $M$ embeds into $Gr_\alpha^{\DD'}$,
which is a subgraph of $Gr^\ast$.
$\eop_{\ref{universal}}$
\end{Proof of the Claim}

We thank Charles Morgan for permitting us to use the following argument 
he showed us:

\begin{Claim}\label{Charles} Let $\DD$ be a normal ultrafilter over $\kappa$
and $A\in \DD$. Suppose that $G$ is $\PPr(\DD)$-generic filter over $V$.
\underline{Then} there is some $G'$ which is $\PPr(\DD)$-generic over $V$
and such that $(\emptyset, A)\in G'$ while $V[G]=V[G']$.
\end{Claim}

\begin{Proof of the Claim} Let $x=x_G=\cup\{s:\,(\exists B\in\DD)(s,B)\in G\}$,
so
\[
G=G_x=\{(s,B)\in \PPr(\DD):\,s\subseteq x_G\subseteq s\cup B\}.
\]
Now we use the Mathias characterisation of Prikry forcing, which says that for
an infinite subset $x$ of $\kappa$ we have that $G_x$ is $\PPr(\DD)$-generic over $V$
iff $x_G\setminus B$ is finite for all $B\in \DD$. Hence $x\setminus A$ is
finite. Let $y=x_G\cap A$, so an infinite subset of $\kappa$ which clearly
satisfies that $y\setminus B$ is finite for all $B\in\DD$. Let $G'=G_y$, so
$G'$ is $\PPr(\DD)$-generic over $V$ and $(\emptyset, A)\in G'$. We have $V[G']
\subseteq V[G]$ because $y\in G$ and $V[G]\subseteq V[G']$ because
$x\setminus y$ is finite.
$\eop_{\ref{Charles}}$
\end{Proof of the Claim}

\begin{Conclusion}\label{homogen} Suppose that $Q$, $\name{\DD}'$, $G$ and $V'$
are as in Claim \ref{universal} and $H$ is a $\PPr(\DD')$-generic filter
over $V'[G]$. \underline{Then} the conclusion of Claim \ref{universal}
holds in $V'[G][H]$.
\end{Conclusion}

\begin{Proof} Since $\EE\deq\{\langle s, B\rangle:\,B\subseteq A^\ast\}$
is in $V'[G]$ and dense in $\PPr(\DD')$, there is $\langle s,B\rangle\in H
\cap \EE$. Let $H^\ast\deq\{\langle t, C\rangle:\,t\supseteq s\}$.
Then $\langle
s, A^\ast\rangle\in H^\ast$ and $H^\ast$ is $\PPr(\DD')/s$-generic over $V'[G]$
with $V'[G][H^\ast]=V'[G][H]$. As $\langle s, A^\ast\rangle$ forces in
$\PPr(\DD')/s$ exactly the same statements as $\langle \emptyset, A^\ast\rangle$
does in $\PPr(\DD')$, the conclusion follows by Claim \ref{universal}.
$\eop_{\ref{homogen}}$
\end{Proof}

\begin{Claim}\label{fits} Suppose that
$\bar{Q}=\langle P_i,\name{Q}_i,\name{A}_i:\,i<\chi
\rangle\in \KK^\ast_{\kappa^+}$ is given by determining
$Q_0$ as in Definition \ref{Q0} and defining $\name{Q}_{\name{\DD}}^i=
\name{Q}_{\name{\DD},\kappa}^{V[G_{P_i}]}$ as defined in Definition \ref{defa},
with $\kappa$ replaced by $\mu$,
and $\name{A}_i=\name{A}^\ast_i$ where $\name{A}^\ast_i$ was  defined
in Claim \ref{properties}(2).

\underline{Then} $\bar{Q}$ is fitted.
\end{Claim}

\begin{Proof of the Claim} We shall take $\name{R}$ to be defined
by Definition \ref{defa}. By Claim
\ref{reflect}, it suffices to give a definition of $\name{f}$ 
satisfying the
requirements of that Claim. Suppose that $P,\name{\DD}$ is such that
$\name{R}[P,\name{\DD}]$,
working in $V^P$ we define
\[
f=f_{[P,\DD]}:\,Q_{\DD}=R_{[P,\DD]}\into\DD
\]
by letting $f(p)\deq B^p$ for $p=(A^p,B^p, u^p, \bar{f}^p)$.
We check that this definition is as required. So suppose that
$\kappa'<\kappa$ is inaccessible and ${\bf g}$ is a $(<\kappa')$-directed
family of conditions in $Q_{\DD}$ with the property that
for all $p\in {\bf g}$ we have $\kappa'\in B^p$. We define $r$
by letting
\[
A^r\deq\bigcup_{p\in {\bf g}}A^p\cup\{\kappa'\},
B^r\deq\bigcap_{p\in {\bf g}}B^p\setminus\{\kappa'\},
u^r\deq\cup_{p\in {\bf g}}u^p,
\]
and for $\alpha\in u^r$, we let $f^r_\alpha\deq\cup_{p\in {\bf g}\,\,\&
\,\,\alpha\in u^p}f^p_\alpha$.
It is easy to check that this condition is as required.
$\eop_{\ref{fits}}$
\end{Proof of the Claim}

\begin{Remark} The inaccessibility of $\kappa'$ was not used in the
Proof of Claim \ref{fits}.
\end{Remark}

\centerline{\em {Proof of the Theorem finished.}}

To finish the proof of the Theorem, in $V_0$ let $\bar{Q}$ be as
in Claim \ref{fits}. By Claim \ref{fits} and the definition
of fittedness,
we can find sequences $\langle p^\ast_i:\,i<\chi\rangle$
and $\langle \alpha_i:\,i<\chi\rangle$
witnessing that $\bar{Q}$ is fitted. Let $\name{\DD}_i\deq
{p^\ast_i(\alpha_i)}$ for $i<\chi$.
If we force
in $V_0$ by
\[
P^\ast
\deq\lim \langle P_{\alpha_i}/(p^\ast_i\rest\alpha_i),\name{Q}_{\DD_i}:\,i<\chi
\rangle,
\]
we obtain a universe $V^\ast$ in which
$\langle \DD_i:\,\cf(i)=\mu^+\rangle$ is an
increasing sequence of normal filters over $\mu$, and $\DD\deq
\bigcup_{i\in S^\chi_{\mu^+}}\DD_i$ is a normal ultrafilter over $\mu$.
For, in $V^{P_{\alpha_i}/(p^\ast_i\rest\alpha_i)}$,
we have that $\DD_i$ is an ultrafilter over $\mu$, and $\cf(\chi)>\mu$,
while the iteration is with $(<\mu)$-supports and $\mu^{<\mu}=\mu$.
Hence every subset if $\mu$ in $V^\ast$ appears as an element of
$V^{P_{\alpha_i}/ (p^\ast_i\rest\alpha_i)}$ for some $i$, and so $\DD$ is
an ultrafilter.

Also, for every $i\in S^\chi_{\mu^+}$ we have that $A_i^\ast\in \DD$.
Let $\name{\DD}$ be a $P^\ast$-name for $\DD$ of $V^\ast$.
Let 
\begin{equation*}
E\deq\left\{\delta<\chi:\,
\begin{split}
(\forall\alpha<\delta)(\exists\beta\in
(\alpha,\delta))[\alpha_\beta=\beta]\mbox{ and }\\
\name{\DD}\cap\PP(\kappa)^{V_0^{P_\beta}}
\mbox{ is a }P_\beta/(p_\beta\rest\beta)\mbox{-name}
\\
\mbox{ and } p_{\beta+1}(\beta)=\DD\cap\PP(\kappa)^{V_0^{P_\beta}}\\
\end{split}
\right\}.
\end{equation*}
Hence $E$ is a club of $\chi$.
Let $\delta\in E\cap S^\chi_{\mu^{++}}$ be larger than $\mu^{+++}$.
Force with $P^\ast\rest\delta$, so obtaining $V_1$ in which $2^{\mu^+}\ge
2^\kappa\ge \mu^{+++}$, as each coordinate of $P^\ast\rest\delta$ adds a
subset of $\mu$, and cardinals are preserved.
In $V_1$ force with the Prikry forcing for $\DD_\delta\deq\bigcup_{i\in
S^\delta_{\mu^+}}\DD_i$.
Let $W\deq V_1[\PPr(\DD_\delta)]$.
For $i\in S^\delta_{\mu^+}$, let ${\rm Gr}^\ast_i$ be a graph obtained in
$W$ satisfying the conditions of Conclusion\ref{homogen} with $\DD_\delta$
in place of
$\DD'$
and $\DD_i$ in place of $\DD$. Let $C$ be a club of $\delta$ of order type
$\mu^{++}$, and let $g$ be its increasing enumeration.

We claim that $W$ is as required, and
that 
\[
\{{\rm Gr}^\ast_{g(i)}:\,i<\mu^{++}\,\,\&\,\,\cf(g(i))=\mu^+\}
\]
are universal for graphs of size $\mu^+$. Clearly the cofinality of $\mu$ in $W$
is
$\aleph_0$ and $\mu$ is a strong limit. 
Suppose that $Gr$ is a graph on $\mu^+$
in $W$ and let $\name{\Gr}$ be a $\PPr(\DD_\delta)$-name for it. Hence, there is
a $i<\mu^{++}$  with $\cf(g(i))=\mu^+$ such that $\name{\Gr}$ is a
$\PPr(\DD_{g(i)})$-name for a graph on $\mu^+$. The conclusion
follows by the choice of ${Gr}^\ast_i$.
$\eop_{\ref{graph}}$ 
\end{Proof}

\begin{Remark} The forcing used in \cite{GiSh 597} also satisfies
the conditions of Claim \ref{reflect}, again with $f(p)=B^p$.
\end{Remark}
\eject

\end{document}